\documentclass[10pt,letter]{article}
\usepackage{amsfonts,epsf,amsmath,amssymb,graphicx,tabularx}

\newtheorem{theorem}{\bf Theorem}[section]
\newtheorem{corollary}[theorem]{\bf Corollary}

\newcommand{\proof}{\noindent{\bf Proof.\ }}
\newcommand{\qed}{\hfill $\blacksquare$ \bigskip}

\begin{document}

\title{\bf{Note on the harmonic index of a graph}}

\author{
Aleksandar Ili\' c \\
Faculty of Sciences and Mathematics, University of Ni\v s, Serbia \\
e-mail: \tt{aleksandari@gmail.com} \\
}

\date{\today}
\maketitle

\begin{abstract}
The harmonic index of a graph $G$ is defined as the sum of weights
$\frac{2}{deg(v) + deg(u)}$ of all edges $uv$ of $E (G)$, where $deg (v)$ denotes the degree of a vertex $v$ in $V (G)$. In this note we generalize results of [L. Zhong, \textit{The harmonic index on graphs}, Appl. Math. Lett. 25 (2012), 561--566] and establish some upper and lower bounds on the harmonic index of $G$.
\end{abstract}

{\bf Key words}: harmonic index, Randi\' c index, trees. \vskip 0.1cm

{{\bf AMS Classifications:} 05C12, 92E10.} \vskip 0.1cm

%%%%%%%%%%%%%%%%%%%%%%%%%%%%%%%%%%%%%%%%%%%%%%%%%%%%%%%
\section{Introduction}
\label{sec:intro}
%%%%%%%%%%%%%%%%%%%%%%%%%%%%%%%%%%%%%%%%%%%%%%%%%%%%%%%

Let $G = (V, E)$ be a connected simple graph with $n = |V|$ vertices
and $m = |E|$ edges. The degree of a vertex $v$ is denoted as $deg
(v)$. 

The Randi\' c index is one of the most successful molecular descriptors in structure-property and structure-activity relationships
studies. The Randi\' c index $R(G)$ is defined as \cite{Ra75}
$$
R(G) = \sum_{uv \in E} \frac{1}{\sqrt{deg (u) deg (v)}}.
$$
The mathematical properties of this graph invariant have been studied extensively (see recent book \cite{GuFu08} and survey \cite{LiSh08}).
Motivated by the success of Randi\' c index, various generalizations and modifications were introduced, 
such as the sum-connectivity index \cite{XiZhTr10,ZhTr09} and the general sum-connectivity index \cite{DuZhTr10,DuZhTr11}.

In this paper, we consider another variant of the Randi\' c index, named the harmonic
index $H(G)$. For a graph $G$, the harmonic index $H(G)$ is defined on the arithmetic mean as
$$
H(G) = \sum_{uv \in E} \frac{2}{deg (u) + deg (v)}.
$$

In \cite{FaMaSa93} the authors considered the relation between the harmonic index and the eigenvalues of graphs. 
Zhong in \cite{Zh11} presented the minimum and maximum values of harmonic index on simple connected graphs and
trees, and characterized the corresponding extremal graphs. In this paper, we generalize some of these results and present alternative simpler proofs of the theorems from \cite{Zh11}.

%%%%%%%%%%%%%%%%%%%%%%%%%%%%%%%%%%%%%%%%%%%%%%%%%%%%%%%
\section{Lower bounds}
\label{sec:lower}
%%%%%%%%%%%%%%%%%%%%%%%%%%%%%%%%%%%%%%%%%%%%%%%%%%%%%%%

Let $G$ be a graph for which holds $deg (u) + deg (v) \leq n$ for all edges $uv \in E (G)$. In triangle-free graphs no two neighboring vertices have a common neighbor, and therefore all triangle-free graphs (trees in particular) belong to this class.

It can be easily seen that
$$
H (G) \geq \sum_{uv \in E (G)} \frac{2}{n} = \frac{2m}{n}.
$$
The equality holds if and only if $deg (u) + deg (v) = n$ for all edges $e = uv \in E (G)$. In particular, for triangle-free graphs the equality holds if and only if $G$ is isomorphic to a complete bipartite graph $K_{a, b}$, $a+b=n$.

\begin{corollary}
Let $T$ be a tree on $n \geq 3$ vertices. Then $H (T) \geq \frac{2 (n - 1)}{n}$ with equality if and only if $T \cong S_n$.  
\end{corollary}

This generalizes Theorem 1 from \cite{Zh11}.
\medskip

Let $\frac{2}{deg (v) + deg (u)}$ be the weight of an edge $e = uv$. Consider an edge $e = uv$ with minimal weight among all edges of $G$, and let $deg (v) = p$ and $deg (u) = q$, with $p, q \geq 2$. Let $x_1, x_2, \ldots, x_p$ be the neighbors of the vertex $v$ and $y_1, y_2, \ldots, y_q$ be the neighbors of $u$. Note that we allow that $x_i \equiv y_j$ for some some vertices $x_i$ and $y_j$. Assume that $G' = G - e$ is a connected graph. Then by using $p + q \geq p + deg (x_i)$ and $p + q \geq q + deg(y_j)$, we have
\begin{eqnarray*}
H (G) - H (G') &=& \left ( \frac{2}{p + q} + \sum_{i = 1}^p \frac{2}{p + deg (x_i)} + \sum_{j = 1}^q \frac{2}{q + deg (y_j)} \right) \\ 
&& - \left ( \sum_{i = 1}^p \frac{2}{p - 1 + deg (x_i)} + \sum_{j = 1}^q \frac{2}{q - 1 + deg (y_j)} \right) \\
&=& \frac{2}{p + q} - \sum_{i = 1}^p \frac{2}{(p +deg(x_i))(p - 1 + deg (x_i))} \\
&& - \sum_{j = 1}^q \frac{2}{(q + deg (y_j)) (q - 1 + deg (y_j))} \\
&\leq& \frac{2}{p + q} - \frac{2 p}{(p + q)(p + q - 1)} - \frac{2 q}{(p + q)(p + q - 1)} \\
&=& \frac{2}{p + q} - \frac{2}{p + q - 1} < 0.
\end{eqnarray*}

Finally, it follows that by removing an edge with the minimal weight, the harmonic index of a graph strictly decreases.
\medskip

%\begin{theorem}\label{sigma} 
%Let $w$ be a vertex of degree $p+1$ in a graph $G$, which is not a
%star, such that $w v_1, w v_2, \ldots, w v_p$ are pendent edges
%incident with $w$ and $u$ is the neighbor of $w$ distinct from $v_1,
%v_2, \ldots, v_p$. We form a graph $G' = \sigma (G, w)$ by removing
%edges $w v_1, w v_2, \ldots, w v_p$ and adding new edges $u v_1, u
%v_2, \ldots, u v_p$. Then
%$$
%H (G) > H (G').
%$$
%\end{theorem}
%
%\proof
%Assume that $deg (w) = x$ and $deg (u) = y > 1$. Then, we have
%$$
%H (G) - H (G') \geq \left( \frac{2}{x + y} + \frac{2 (x- 1)}{x + 1} \right) - \frac{2x}{x + y} > 0,
%$$
%since the weights of the other edges incident to $u$ in $G'$ decreased.
%\qed

The first Zagreb index is one of the oldest molecular
structure-descriptors \cite{GuTr72,GuDa04} and defined as the sum
of squares of the degrees of the vertices
$$
M_1 (G) = \sum_{v \in V} deg (v)^2.
$$

By using Cauchy-Schwarz inequality, we get
$$
\left( \sum_{uv \in E (G)} \frac{2}{deg (u) + deg (v)} \right) \left( \sum_{uv \in E (G)} \frac{deg (u) + deg (v)}{2} \right)  \geq \left( \sum_{e \in E (G)} 1 \right)^2,
$$
or equivalently 
$$
H (G) \geq \frac{2m^2}{M_1 (G)},
$$
with equality if and only if $deg (u) + deg (v)$ is a constant for each edge $e = uv$ from $G$.

\newpage

%%%%%%%%%%%%%%%%%%%%%%%%%%%%%%%%%%%%%%%%%%%%%%%%%%%%%%%
\section{Upper bounds}
\label{sec:upper}
%%%%%%%%%%%%%%%%%%%%%%%%%%%%%%%%%%%%%%%%%%%%%%%%%%%%%%%

\begin{theorem}
\label{thm-pi} Let $w$ be a vertex of a nontrivial connected graph $G$. For nonnegative integers
$p$ and $q$, let $G (p, q)$ denote the graph obtained from $G$ by attaching to the vertex $w$ pendent
paths $P = w v_1 v_2 \ldots v_p$ and $Q = w u_1 u_2 \dots u_q$ of lengths $p$ and $q$,
respectively. Then
$$
H (G (p, q)) < H (G (p + q, 0)).
$$
\end{theorem}

\proof
Without loss of generality, assume that $p \geq q > 0$. Since $G$ is a nontrivial connected graph, we have $x = deg (w) > 2$, and after transformation the vertex degree of $w$ decreases by one and the weights of the edges in $G$ either remain the same or decrease (the later are the edges adjacent to $w$). We will consider the difference $\Delta = H (G (p + q, 0)) - H (G (p, q))$ in three cases.

{\bf Case 1. } $p = q = 1$ 
$$
\Delta > \frac{1}{x + 1} + \frac{1}{3} - \frac{1}{x + 1} - \frac{1}{x + 1} = \frac{1}{3} - \frac{1}{x + 1} > 0.
$$

{\bf Case 2. } $p > q = 1$ 
$$
\Delta > \frac{1}{x + 1} + \frac{1}{3} + \frac{p - 1}{4} - \frac{1}{x + 1} - \frac{1}{x + 2} - \frac{p - 2}{4} - \frac{1}{3} = \frac{1}{4} - \frac{1}{x + 2} > 0.
$$

{\bf Case 3. } $p \geq q > 1$ 
$$
\Delta > \frac{1}{x + 1} + \frac{1}{3} + \frac{p + q - 2}{4} - \frac{2}{x + 2} - \frac{p + q - 4}{4} - \frac{2}{3} = \frac{1}{6} - \frac{2}{x + 2} + \frac{1}{x + 1} > 0.
$$

This completes the proof. \qed

By repetitive application of this transformation on branching vertices that are on the largest distance from the center of $T$, we have
\begin{corollary}
Let $T$ be a tree on $n \geq 3$ vertices. Then $H (T) \leq \frac{4}{3} + \frac{n-3}{4}$ with equality if and only if $T \cong P_n$.  
\end{corollary}

This generalizes Theorem 2 from \cite{Zh11}. \medskip

For $n \geq 7$, using Theorem \ref{thm-pi}, one can conclude that the tree with the second maximum value of harmonic index is of the form  $T (a, b, c)$, composed from the branching vertex $v$ of degree 3 with three pendent paths $P_a, P_b, P_c$ attached to $v$ ($a + b + c + 1 = n$). Without loss of generality, assume that $a \geq b \geq c$ and

\begin{itemize}
\item $b = c = 1$: $H (T(a, b, c)) = \frac{1}{5} + \frac{2}{4} + \frac{1}{3} + \frac{n - 5}{4}$
\item $b > c = 1$: $H (T(a, b, c)) = \frac{2}{5} + \frac{2}{3} + \frac{1}{4} + \frac{n - 6}{4}$
\item $b \geq c > 1$: $H (T(a, b, c)) = \frac{3}{5} + \frac{3}{3} + \frac{n - 7}{4}$
\end{itemize}

By direct verification, we get that trees $T (a, b, c)$ with $a \geq b \geq c > 1$ have the second maximum value of the harmonic index among trees on $n$ vertices.

\bigskip {\bf Acknowledgment. } This work was supported by Research
Grants 174010 and 174033 of Serbian Ministry of Education and Science.

\end{document}